\documentclass[sn-mathphys-num]{sn-jnl}

\usepackage{upgreek}
\usepackage{amsmath,amssymb,amsthm,amsfonts,mathrsfs}
\usepackage{caption}
\usepackage{subcaption}
\usepackage{xcolor}
\usepackage{ dsfont }
\usepackage{ stmaryrd }
\usepackage{pifont}
\usepackage{ textcomp }
\usepackage{tabularx}
\usepackage{multirow}

\newtheorem{Ex}{Example}[section]
\usepackage{graphicx}
\usepackage{amsmath,amssymb,amsfonts}%
\usepackage{mathrsfs}%
\usepackage[title]{appendix}%
\usepackage{xcolor}%
\usepackage{textcomp}%
\usepackage{manyfoot}%
\usepackage{booktabs}%
\usepackage{algorithm}%
\usepackage{algorithmicx}%
\usepackage{algpseudocode}%
\usepackage{listings}%
\usepackage{subcaption} 
\definecolor{myNewColorA}{RGB}{126,12,110}
\definecolor{myNewColorB}{RGB}{165,85,154}
\definecolor{myNewColorC}{RGB}{203,158,197}
\usepackage{hyperref}
\usepackage{natbib}

\begin{document}
	
	\title[Fractional Order Sunflower Equation: Stability, Bifurcation and Chaos]{Fractional Order Sunflower Equation: Stability, Bifurcation and Chaos}
	
	
	\author[1]{\fnm{Deepa} \sur{Gupta}}\email{deepanew96@gmail.com}
	
	\author*[2]{\fnm{Sachin} \sur{Bhalekar}}\email{sachinbhalekar@uohyd.ac.in}

	\affil[1,2]{\orgdiv{School of Mathematics and Statistics}, \orgname{University of Hyderabad}, \orgaddress{\city{Hyderabad}, \postcode{500046}, \state{Telangana}, \country{India}}}

	
	\abstract	{The sunflower equation describes the motion of the tip of a plant due to the auxin transportation under the influence of gravity. This work proposes the fractional-order generalization to this delay differential equation. The equation contains two fractional orders and infinitely many equilibrium points. The problem is important because the coefficients in the linearized equation near the equilibrium points are delay-dependent. We provide a detailed stability analysis of each equilibrium point using linearized stability. We find the boundary of the stable region by setting the purely imaginary value to the characteristic root. This gives the conditions for the existence of the critical values of the delay at which the stability properties change. We observed the following bifurcation phenomena: stable for all the delay values, a single stable region in the delayed interval, and a stability switch. We also observed a multi-scroll chaotic attractor for some values of the parameters.}

	
	
	
	\maketitle
	
	\section{Introduction}
	Mathematical analysis is the fundamental tool used to study complex Biological systems \cite{edelstein2005mathematical,segel2013primer}. These systems can be modeled using differential or difference equations. If the ``time" variable in these systems is continuous, then one can use the ordinary or partial differential equation \cite{taubes2008modeling,logan2004partial}. Such equations may be improved by including ``time-delay," which gives a better fit for the ``system-memory" \cite{smith2011introduction,rihan2021delay}. Control and bifurcation studies in the delayed predator-prey system are carried out by Ou et al. in \cite{ou2024hopf}. Karkri et al. \cite{el2024stability} analyzed the delayed equation arising in the infection dynamics. Hopf bifurcation analysis of the neural networks is presented in \cite{li2023exploring}. Various applications of delay differential equations (DDE) in biological systems are available in a special issue edited by Rihan et al. \cite{rihan2018applications}.
	\par The obvious nonlocality of the Biological systems can be captured in the model with the help of ``fractional order derivative (FOD)" \cite{magin2004fractional,bhalekar2011fractional,arshad2019fractional}. If the order of the derivative is a non-integer (e.g., a positive real number or a complex number with a positive real part), then it is termed as FOD \cite{podlubny1998fractional}. Mathematicians are working on derivatives whose order depends on time or is distributed over some interval \cite{sun2019review,mainardi2008time}. Magin \cite{magin2010fractional} used fractional derivatives to study the complex models in biological systems. The stability of the fractional order delay differential equations (FDDE) is studied by Bhalekar and coworkers \cite{bhalekar2016stability,bhalekar2013stability,bhalekar2019analysis,bhalekar2022stability}. Latha et al.  \cite{latha2017fractional} proposed an Ebola infection model using FDDEs. The optimal control of cancer treatment using the FDDE model can be found in \cite{sweilam2020fractional}. FDDEs are used to study the population dynamics in \cite{xu2023exploring}. Fractional order three-triangle multi-delayed neural networks are discussed by Xu et al. \cite{xu2023bifurcation}. Complex valued neural networks involving FDDEs are presented by Rakkiyappan et al. \cite{rakkiyappan2017stability}. Fractional order octonion-valued neural networks with delay are investigated by  Udhayakumar and Rakkiyappan \cite{kandasamy2020hopf}. Bhalekar et al. \cite{bhalekar2011fractional} proposed the fractional order bloch equation with delay.
	
	\par In 1967, Israelsson and Johnsson \cite{israelsson1967theory} proposed a model known as the ``sunflower equation" explaining the helical movements (circumnutations) of the apex of sunflower plants. The theory concerns the interplay between gravity and the growth hormone (auxin). The time delay arises because the hormone takes some time to spread in the plant body (geotropic reaction time for the hypocotyls). Somolinos \cite{somolinos1978periodic} in 1978 carried out the rigorous mathematical analysis of the sunflower equation. Oscillations in this equation are studied by Kulenovi{\'c} and Ladas \cite{kulenovic1988oscillations}.
	\par In this work, we generalize the sunflower equation to include the fractional order derivatives. The fractional order gives the flexibility to select the order $\alpha$ in the interval $(0, 1]$ in contrast with the fixed value $\alpha=1$ in the classical model. This flexibility is extremely useful in fitting the experimental data with the mathematical model. Furthermore, unlike the classical integer-order derivatives, the FDO is a nonlocal operator. This nonlocality is helpful in modeling the memory properties of the natural systems. Thus, our model is an improved version of the sunflower equation described in the literature. The system possesses infinitely many equilibrium points. We provide the stability analysis of all these equilibrium points and discuss the possible types of bifurcations in detail which is not available in the literature discussed above. Furthermore, we present the chaotic solutions of this system. The novelty of this work is that it investigates a multi-scroll chaotic attractor in the proposed model.
	
	\par The rest of the paper is organized as follows: Section \eqref{sec4.1} details the sunflower equation and its fractional-order counterpart. In Section \eqref{sec4.4}, we describe the sunflower equation's stability and bifurcation analysis of equilibrium points.
	Chaos in the proposed model is studied in Section \eqref{sec4.5}.  Validation of results is done in Section \eqref{sec4.6}. Section \eqref{sec4.7} presents the conclusions. 
	
	\section{The Sunflower Equation}\label{sec4.1}
	The sunflower equation \cite{israelsson1967theory} described by (\ref{sun})  is a modeling nonlinear delay differential equation that defines the helical movement of the tip of a growing plant that accumulates growth hormone (auxin).
	
	\begin{equation} 
		\dfrac{\tau}{l}\ddot{x}(t)+\dot{x}(t)=\dfrac{-m}{l}\sin(x(t-\tau)) \label{sun}	
	\end{equation}
	
	where $m$, $l$ and delay $\tau$ are positive numbers. Now, we generalize it to the fractional order case as:
	\begin{equation}\label{eq4.1}
		\dfrac{\tau}{l}D^{2\alpha}x(t)+D^\alpha x(t)=\dfrac{-m}{l}\sin(x(t-\tau)) ,  \quad 0<\alpha\leq 1.	
	\end{equation}
	Here, $D^\alpha$ and $D^{2\alpha}$ represent Caputo fractional derivatives \cite{bhalekar2013stability,diethelm2002analysis,kilbas2006theory,lakshmanan2011dynamics,podlubny1998fractional}.

	Note that $x^*$ is an equilibrium point of equation \eqref{eq4.1} if and only if $\sin(x^*)=0$ as we have $D^\alpha x^*=0$ and $D^{2\alpha} x^*=0$.
	So, equation \eqref{eq4.1} has infinitely many equilibrium points given by $x_{1,n}^{*}=2 n \pi$ and $x_{2,n}^{*}=(2n+1)\pi$, $n\in \mathbb{Z}$.

	By taking a small perturbation near the equilibrium point and using Taylor's approximations, we get the local linearization of equation \eqref{eq4.1} as \begin{equation}\label{eq4.2}
		\dfrac{\tau}{l}D^{2\alpha}x(t)+D^\alpha x(t)=\dfrac{-m}{l}x(t-\tau),
	\end{equation}
	near the equilibrium points $x_{1,n}^{*}$ and	
	\begin{equation}\label{eq4.3}
		\dfrac{\tau}{l}D^{2\alpha}x(t)+D^\alpha x(t)=\dfrac{m}{l}x(t-\tau),
	\end{equation}
	near the equilibrium points $x_{2,n}^{*}$.	
	\section{Stability Analysis}\label{sec4.4}	
	If we consider the non-delayed equation \eqref{eq4.2} then we get $D^\alpha x(t)=\dfrac{-m}{l}x(t)$ which implies that the equilibrium points $x_{1,n}^*$ are stable at $\tau=0$ \cite{matignon1996stability}. Similarly, the equilibrium points $x_{2,n}^*$ are unstable at $\tau=0$ as the equation \eqref{eq4.3} gets reduced to $D^\alpha x(t)=\dfrac{m}{l}x(t)$. \\
	\subsection{Stability results for the equilibrium points $x_{1,n}^*=2n\pi$}
\begin{figure}
	\subfloat[$x_{1,n}^*$ is stable for all $\tau$ when there is no intersection between $g(\tau)$ and $\tau$]{%
		\includegraphics[scale=0.7]{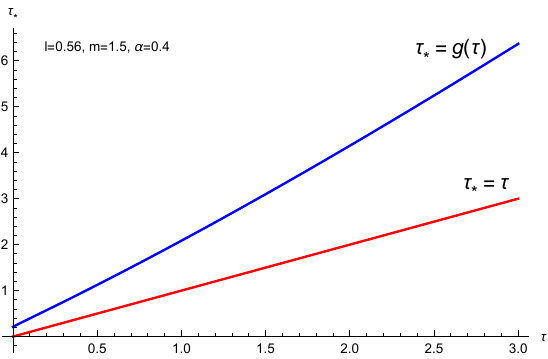}
	}\hspace{0.1cm}
	\subfloat[Bifurcation of the stable region from the stability switch]{%
		\includegraphics[scale=0.75]{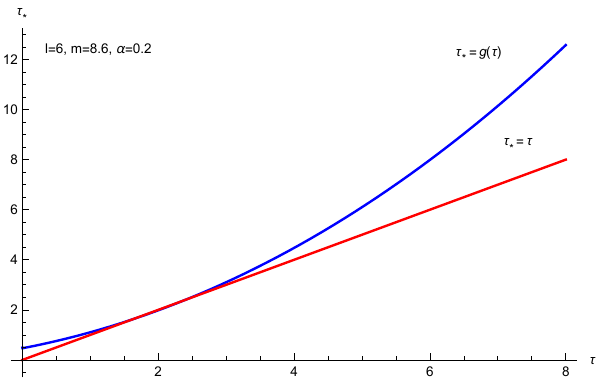}
	}\hspace{0.1cm}
	\subfloat[Two intersections between $g(\tau)$ and $\tau$ results the stability switch S-U-S]{%
		\includegraphics[scale=0.8]{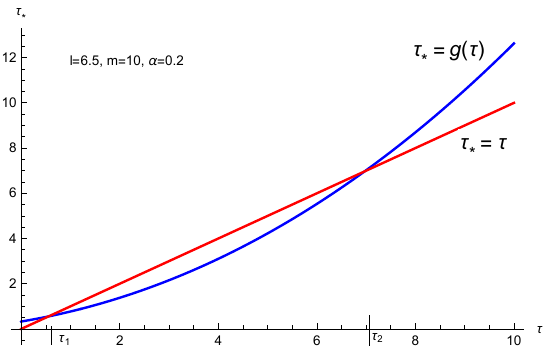}
	}\hspace{0.1cm}
	\subfloat[Only one intersection between $g(\tau)$ and $\tau$ gives the single stable region]{%
		\includegraphics[scale=0.8]{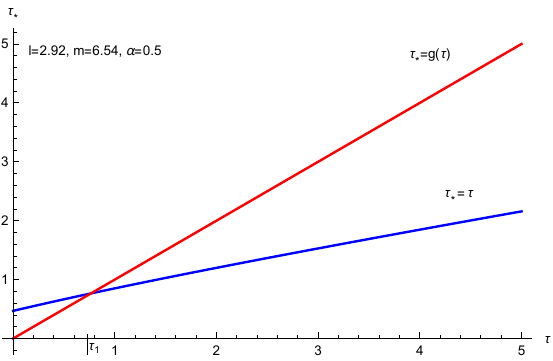}
	}
	\caption{ Different behaviors of the critical curve $\tau_*=g(\tau)$ in the $\tau\tau_*$ plane}
	\label{figure4.3}
\end{figure}

	Now, let us consider the equilibrium points $x_{1,n}^*$ and $\tau>0$. \\
	By using Laplace transform, the characteristic equation of \eqref{eq4.2} is: 
	\begin{equation}\label{eq4.4}
		\dfrac{\tau}{l}	\lambda^{2\alpha}+\lambda^\alpha+\dfrac{m}{l}\exp(-\lambda\tau)=0.
	\end{equation}
	We have a change in stability only when the root $\lambda=u+iv$ of equation \eqref{eq4.4} crosses the imaginary axis. \\
	Therefore, by putting $\lambda=iv$, $v>0$ in  the equation \eqref{eq4.4}, we get the boundary of the stable region
	i.e. $	\dfrac{\tau}{l}	(i v)^{2\alpha}+(i v)^\alpha+\dfrac{m}{l}\exp(-iv\tau)=0.$ \\
	Separating the real and imaginary parts, we get,
	\begin{equation}\label{eq4.5}
		\dfrac{\tau}{l}v^{2\alpha}\cos(\alpha\pi)+v^{\alpha}\cos\Big(\dfrac{\alpha\pi}{2}\Big)=\dfrac{-m}{l}\cos(v\tau)
	\end{equation}
	and\begin{equation}\label{eq4.6}
		\dfrac{\tau}{l}v^{2\alpha}\sin(\alpha\pi)+v^{\alpha}\sin\Big(\dfrac{\alpha\pi}{2}\Big)=\dfrac{m}{l}	\sin(v\tau).
	\end{equation}
	Now, by squaring and adding equations \eqref{eq4.5} and \eqref{eq4.6}, we get  
	\begin{equation}\label{eq4.7}
		l^2v^{2\alpha}+\tau^2v^{4\alpha}+2 l v^{3\alpha}\tau\cos(\dfrac{\alpha\pi}{2})-m^2=0
	\end{equation}
	Since $l$ and $m$ are positive numbers, we get only one positive root $v^\alpha$ of equation \eqref{eq4.7} given in the Data Set 1 available at \url{https://drive.google.com/drive/folders/1jOuemmKoSxZfzFSRlotf94YYp5nJI-iy?usp=sharing}. By putting this value of $v$ in the equation \eqref{eq4.5}, we get a critical value of delay $\tau_*$ where we have a change in stability for the equilibrium point $x_{1,n}^*$.\\
	The existence of positive root $v^\alpha$ implies the existence of critical value of delay $\tau_*$. If $Re\Big[\dfrac{d\lambda}{d\tau}\Big|_{u=0}\Big]<0$ at $\tau_*$ then $\exists$ some characteristic root on the right half complex plane moving from the right half complex plane to the left half plane. Since the equilibrium $x_{1,n}^*$ is stable at $\tau=0$, there does not exist any root on the right half plane. This contradiction implies that $Re\Big[\dfrac{d\lambda}{d\tau}\Big|_{u=0}\Big]>0$ at $\tau_*$. Therefore, if $\tau>\tau_*$ then equilibrium point $x_{1,n}^*$ is unstable. Also, there exist only one positive root of $v^\alpha$ where $Re\Big[\dfrac{d\lambda}{d\tau}\Big|_{u=0}\Big]>0$ and the lowest critical value of $\tau$ where we have change in stability is given by $\tau_*$. See Section (6) in \cite{bhalekar2024stability} for more details. Since the coefficient in equation \eqref{eq4.4} depends on $\tau$, the expression for $\tau_*$ also depends on $\tau$, say $\tau_*=g(\tau)$ and is given in the above Data Set 1.\\
	Now, note that if the curve $g(\tau)$ does not meet the line $\tau_*=\tau$ in the $\tau\tau_*$plane then $x_{1,n}^*$ is asymptotically stable (cf. Figure \eqref{figure4.3}(a)). If it meets twice, then $x_{1,n}^*$ will generate a stability switch (SS) as shown in Figure \eqref{figure4.3}(c). So, there exists $\tau_1$ and $\tau_2$ where $g(\tau_1)=\tau_1$ and $g(\tau_2)=\tau_2$ such that when $\tau\in[0,\tau_1)$ then $x_{1,n}^*$ is asymptotically stable, if $\tau\in(\tau_1,\tau_2)$ then $x_{1,n}^*$ is unstable and if $\tau>\tau_2$ again we get $x_{1,n}^*$ asymptotically stable (cf. Figure \eqref{figure4.3}(c)). If $g(\tau)$ cuts only once the line $\tau=\tau_*$, then we have a single stable region (SSR) as given in Figure \eqref{figure4.3}(d).
	So, if $\tau_1$ is the only intersection point where $g(\tau_1)=\tau_1$ then  then $x_{1,n}^*$ is asymptotically stable when $\tau<g(\tau_1)$ whereas $\tau>g(\tau_1)$ implies that $x_{1,n}^*$ is unstable.
	
	\subsection{Bifurcation analysis of $x_{1,n}^*$}\label{subsec4.2}
\begin{figure}
	\subfloat[Bifurcation curve in the $lm-$plane for $\alpha=0.1$, $\alpha=0.2$ and $\alpha=0.3$ (not to the scale)]{%
		\includegraphics[scale=0.7]{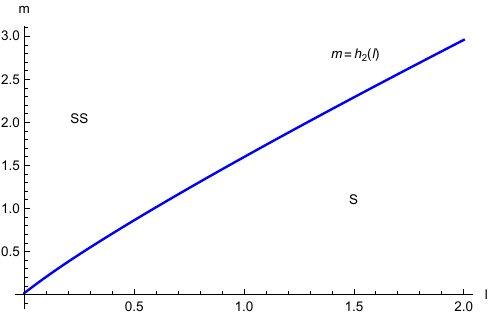}
	}\hspace{0.1cm}
	\subfloat[Bifurcation curve in the $lm-$plane for $\alpha=0.4$]{%
		\includegraphics[scale=0.7]{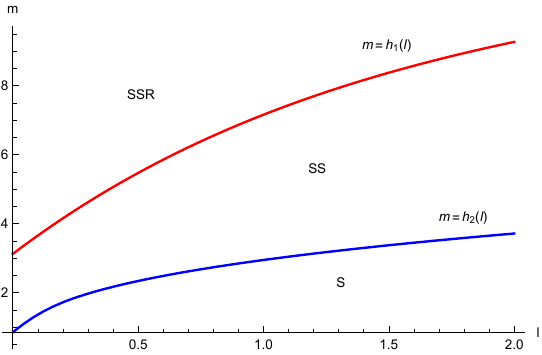}
	}
	\caption{Bifurcation curves in the $lm-$plane for different values of $\alpha$ separating the stable region (S), the stability switch (SS) and the single stable region (SSR)}
	\label{fig4.4}
\end{figure}

	\begin{itemize}
		\item  For $\alpha=0.1$, $0.2$ and $0.3$, we observed only two behaviors of $x_{1,n}^*$ viz. stable (Figure \eqref{figure4.3}(a)) and stability switch (SS) (Figure \eqref{figure4.3}(c)). At the bifurcation of these two behaviors, the curve $\tau_*=\tau$ becomes tangent to the curve $\tau_*=g(\tau)$ (cf. Figure \eqref{figure4.3}(b)). The values of $l$ and $m$ (cf. Figure \eqref{fig4.4}) at such tangent form a curve $m=h_2(l)$ in $lm-$plane that bifurcates stable region with the stability switch. 	\\
		
		\item For $\alpha=0.4$, we get all the stability behaviors viz. S, SS and SSR.\\
		
		If we take $m>h_2(l)$, then there are two intersections between $g(\tau)$ and $\tau$, which results in the SS region. The second intersection point goes away from the first as we increase $m$ further. At the another bifurcation $m=h_1(l)$ (see Figure \eqref{fig4.4}(d)), the second intersection point $\rightarrow\infty$ and vanishes. \\
		Therefore, for $m>h_1(l)$, there is only one intersection between $g(\tau)$ and $\tau$, and we get the SSR. This gives another bifurcation curve $m=h_1(l)$ separating SS from SSR in the $lm-$plane.\\ 
		
		\item For  $1/2\leq\alpha<1$, we observed that the curve $\tau_*=g(\tau)$ and $\tau_*=\tau$ have only one intersection (Figure \eqref{figure4.3}(d)). Thus, there is SSR for $x_{1,n}^*$ and no bifurcation is observed.
	\end{itemize}

	\subsection{Stability results for the equilibrium points $x_{2,n}^{*}=(2n+1)\pi$}\label{sec4.3}
	The characteristic equation corresponding to  $x_{2,n}^{*}$ is 
	\begin{equation}\label{another characteristic equation}
		\frac{\tau}{l} \lambda^{2\alpha}+\lambda^\alpha-\frac{m}{l}\exp(-\lambda\tau)=0.
	\end{equation} 
	If we take $P(\lambda,\tau)=\frac{\tau}{m}\lambda^{2\alpha}+\frac{l}{m}\lambda^\alpha$ and $Q(\lambda,\tau)=\exp(-\lambda\tau)$ then the characteristic equation \eqref{another characteristic equation} can be rewritten as $P(\lambda,\tau)=Q(\lambda,\tau)$.\\
	So, when the graphs of $P(\lambda,\tau)$ and $Q(\lambda,\tau)$ meet at $\lambda_0$ then we get a characteristic root $\lambda_0$ of \eqref{another characteristic equation}, for some $\tau>0$.\\
	Now, $0<\exp(-\lambda\tau)\leq1$ for any real $\lambda>0$ and $\tau\geq0$. So, the range of $Q(\lambda,\tau)=(0,1]$. If we could show that the interval $(0,1]$ is also contained in the range of $P(\lambda,\tau)$ for $\lambda>0$  and $\tau>0$, then there always exists a positive real root $\lambda$ of the characteristic equation \eqref{another characteristic equation}. Since, we have $P(0,\tau)=0$ and $\frac{dP}{d\lambda}=\frac{\tau}{m}(2\alpha)\lambda^{2\alpha-1}+\frac{l}{m}\alpha \lambda^{\alpha-1}>0$, $\lim_{\lambda\rightarrow\infty}P(\lambda,\tau)\rightarrow\infty$. So, the graph of $P(\lambda,\tau)$ is monotonically increasing and contains the interval $(0,1]$. So, we always have a real positive root $\lambda$ of the characteristic equation \eqref{another characteristic equation} for any $\tau>0$. Hence, the equilibrium point $x_{2,n}^{*}$ is always unstable for all $\tau\geq0$.
	
	\section{Chaos In The Sunflower Equation} \label{sec4.5}
\begin{figure}
		\subfloat[Stable solution for $\tau=4$]{%
			\includegraphics[scale=0.4]{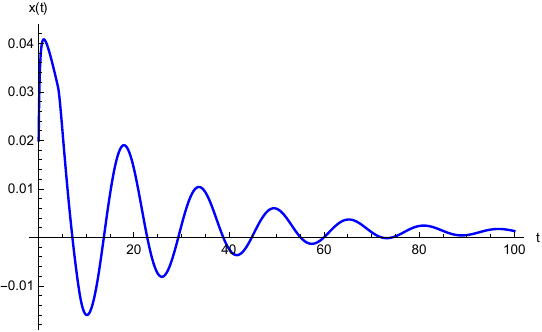}
		}\hspace{0.1cm}
		\subfloat[Unstable solution for $\tau=6$]{%
			\includegraphics[scale=0.4]{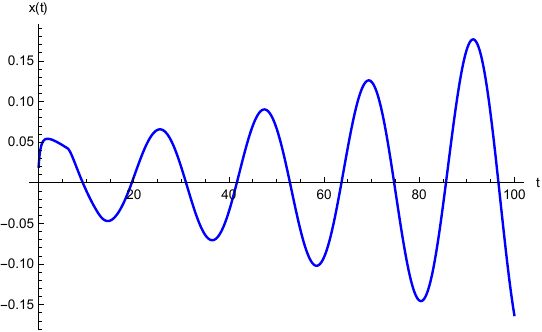}
		}\hspace{0.1cm}
		\subfloat[Asymptotic one-cycle for $\tau=8$]{%
			\includegraphics[scale=0.4]{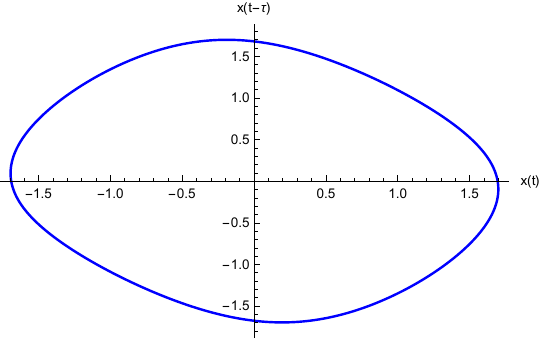}
		}\hspace{0.1cm}
		\subfloat[Asymptotic two-cycle for $\tau=14$]{%
			\includegraphics[scale=0.4]{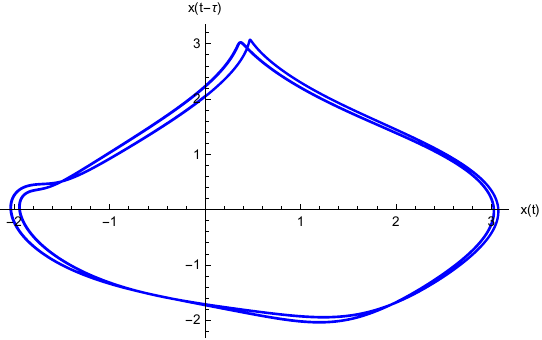}
		}
		\hspace{0.1cm}
		\subfloat[Asymptotic four-cycle for $\tau=15$]{%
			\includegraphics[scale=0.4]{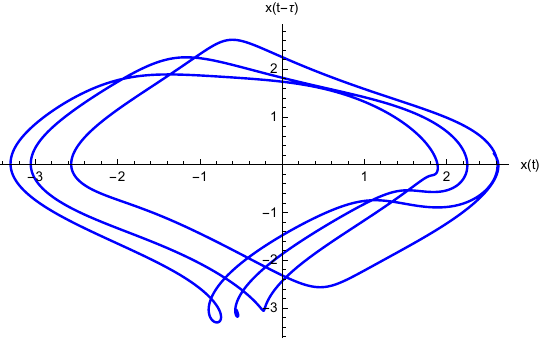}
		}
		\hspace{0.1cm}
		\subfloat[Chaotic attractor for $\tau=20$]{%
			\includegraphics[scale=0.4]{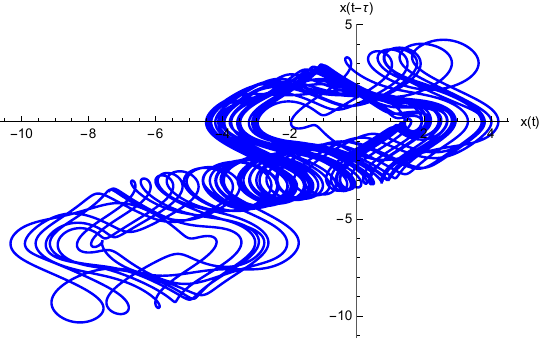}
		}
		\caption{Period doubling route to chaos for $l=14$, $m=5.6$ and $\alpha=0.85$}
		\label{fig4.5}
	\end{figure}

\begin{figure}
		\centering
		\includegraphics{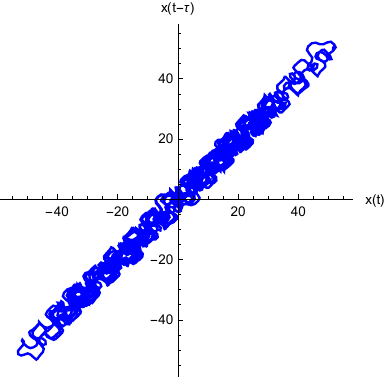}
		\caption{Infinite scroll chaotic attractor }
		\label{fig4.8}
	\end{figure}

	If we take $l=14$, $m=5.6$ and $\alpha=0.85$ then we get the critical value  $\tau_1=5.16433$ (Figure \eqref{figure4.3}(d)) where we have $g(\tau)=\tau$. So, if we take $\tau<\tau_1$ we get stable solution near $x_{1,n}^*$. The stable solution for $\tau=4$ with initial data $x(t)=6.9$ and $\dot{x}(t)=2.5$, $-\tau< t\leq0$  is shown in Figure \eqref{fig4.5}(a). We get the unstable solution for $\tau>\tau_1$ as shown in Figure \eqref{fig4.5}(b) with $\tau=6$.\\
	We used the predictor-corrector method provided in \cite{bhalekar2019analysis,daftardar2014new} for the numerical simulations in this work.\\
	If we further increase the delay $\tau$, e.g., $\tau=8$, we get an asymptotic periodic solution converging to one-cycle (cf. Figure \eqref{fig4.5}(c)). For $\tau=14$ and $\tau=15$, we get asymptotic two-cycle and four-cycle, respectively, as shown in Figures \eqref{fig4.5}(d) and \eqref{fig4.5}(e). The period doubling leads to chaos \eqref{fig4.5}(f). The maximum Lyapunov exponent (MLE) for this fractional-order case with delay $\tau=20$ is $0.2383$. The positive value of the MLE confirms the chaotic oscillations. We used the program by Kodba, Perc, and Marhl \cite{kodba2004detecting}, which is based on Wolf's algorithm \cite{wolf1986quantifying} to find the MLE.  Figure \eqref{fig4.8} shows the infinite-scroll chaotic attractor for $\alpha=1$, $l=14$, $m=5.6$ and $\tau=20$. The MLE, in this case, is $0.3539$.\\
	Note that the fractional order systems do not have exactly periodic orbits \cite{kaslik2012non}. However, we can have asymptotic-periodic orbits or limit cycles, as observed in this work. 
	\section{Examples}\label{sec4.6}
\begin{figure}
		\subfloat[Stable solution for $\tau=4$ and $m=1$]{%
			\includegraphics[scale=0.3]{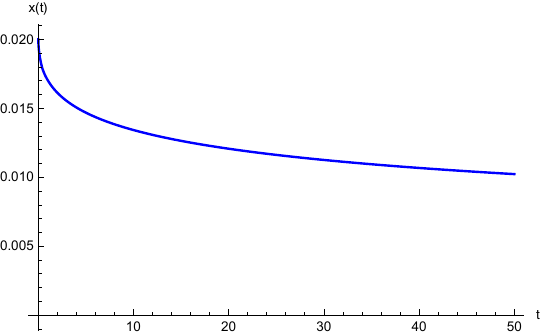}
		}\hspace{0.1cm}
		\subfloat[Stable solution for $\tau=0.4$ and $m=6$]{%
			\includegraphics[scale=0.3]{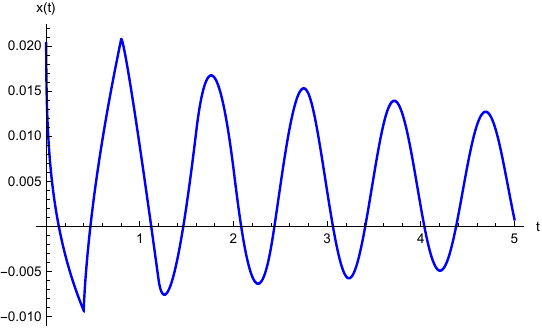}
		}\hspace{0.1cm}
		\subfloat[Unstable solution for $\tau=0.7$ and $m=6$]{%
			\includegraphics[scale=0.3]{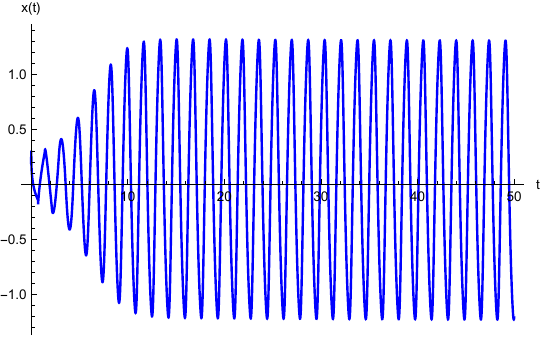}
		}\hspace{0.1cm}
		\subfloat[Stable solution for $\tau=12$ and $m=6$]{%
			\includegraphics[scale=0.3]{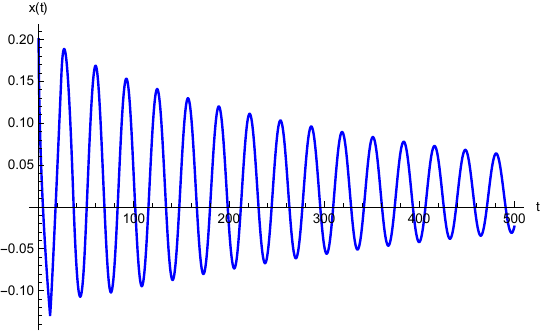}
		}
		\caption{Figures of Example \eqref{eg4.1}}
		\label{fig4.6}
	\end{figure}

	\begin{Ex}\label{eg4.1}
		Figure  \eqref{fig4.4} shows that there is only one bifurcation curve in the $lm-$plane separating the region S from SS for $\alpha=0.1$, $0.2$, and $0.3$.
		
		\begin{itemize}
			\item	So, if we take $\alpha=0.3$ and $l=3$ with the initial data $x(t)=0.02$ for $t\in(-\tau,0]$ near $x_{1,0}^*=0$  then along the bifurcation curve we get the critical value $m_*=h_2(l)=5.3092$. 
			\begin{itemize}
				\item	If we take $m=1<m_*$ then the equilibrium point $0$ is stable $\forall\tau\geq0$. Figure \eqref{fig4.6}(a) shows the stable solution for $\tau=4$.
				\item	Now, if we take $m=6>m_*$ then we are in the stability switch region. The two critical values of $\tau$ are $\tau_1=0.567501$ and $\tau_2=10.133$ (see Figure \eqref{figure4.3}(c)) where $g(\tau)=\tau$. If we take $\tau<\tau_1$, we get stable solution near $0$ (cf. Figure \eqref{fig4.6}(b) for $\tau=0.4$), if we take $\tau_1<\tau<\tau_2$ we get unstable solution near $0$. So, if we take $\tau=0.7$ and all the other parameters are fixed, then one of the characteristic roots with positive real parts is $0.0425373+3.65101 i$, which is sufficient for the instability. Now, if we take $\tau>\tau_2$ (cf. Figure \eqref{fig4.6}(d) for $\tau=12$) again we get stable solution near $0$.
			\end{itemize}
		\end{itemize}
	\end{Ex}
	
\begin{figure}
		\subfloat[Stable solution for $\tau=0.08$ and $m=1$]{%
			\includegraphics[scale=0.4]{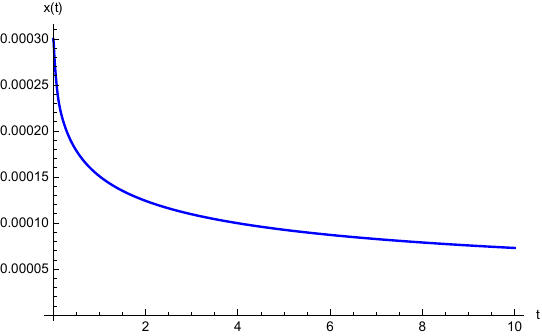}
		}\hspace{0.1cm}
		\subfloat[Stable solution for $\tau=0.4$ and $m=3.2$]{%
			\includegraphics[scale=0.4]{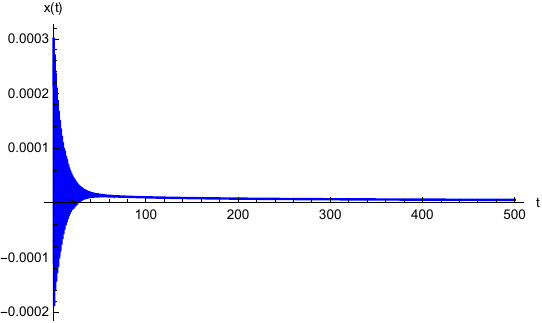}
		}\hspace{0.1cm}
		\subfloat[Unstable solution for $\tau=0.8$ and $m=3.2$ ]{%
			\includegraphics[scale=0.4]{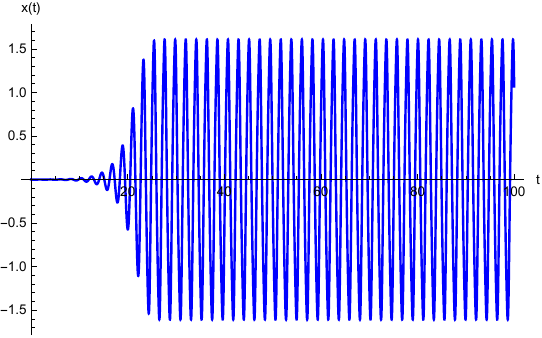}
		}\hspace{0.1cm}
		\subfloat[Stable solution for $\tau=12$ and $m=3.2$]{%
			\includegraphics[scale=0.4]{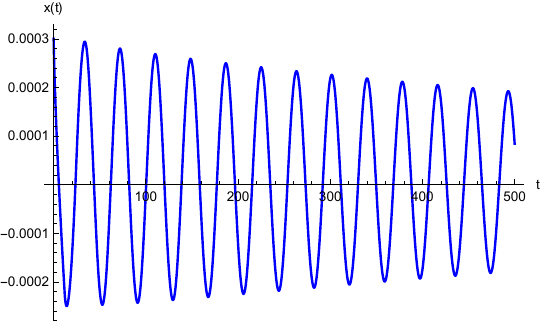}
		}
		\hspace{0.1cm}
		\subfloat[Stable solution for $\tau=0.01$ and $m=8$]{%
			\includegraphics[scale=0.4]{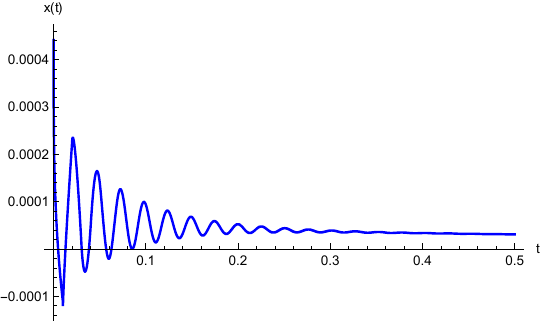}
		}\hspace{0.1cm}
		\subfloat[Unstable solution for $\tau=0.03$ and $m=8$]{%
			\includegraphics[scale=0.4]{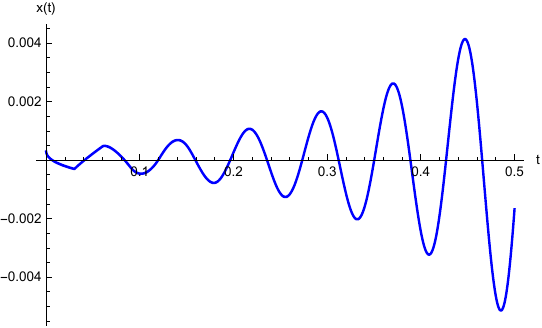}
		}
		\caption{Figures of Example \eqref{eg4.2} with $l=1$ and $\alpha=0.4$}
		\label{fig4.7}
	\end{figure}

	\begin{Ex}\label{eg4.2}
		Figure \eqref{fig4.4}(d) shows that there are all three types of behaviors viz. S, SS, and SSR for $\alpha=0.4$ in a neighborhood of the equilibrium points $x_{1,n}^*$.
		\begin{itemize}	
			\item So, if we fix $l=1$ then we get the critical values $m_{1*}=2.95108$ and $m_{2*}=7.16$ along the curves $m=h_2(l)$ and $m=h_1(l)$, respectively.  We take initial data as $x(t)=0.0003$ $\forall t\in(-\tau,0]$ near  $x_{1,0}^*=0$.
			\begin{itemize}
				\item If we take $m=1<m_{1*}$, we get a stable solution near $0$  for all $\tau>0$. Figure \eqref{fig4.7}(a) shows stable solution near $0$ for $\tau=0.08$ .
				\item Now, if we take $m=3.2 \in(m_{1*}, m_{2*})$, then it is in the stability switch region, and we get $\tau_1=0.616608$ and $\tau_2=10.733$ given in the Figure \eqref{figure4.3}(c). So, we get stable solution for $\tau\in[0, \tau_1)$, unstable solution for $\tau\in(\tau_1, \tau_2)$ and again stable solution for $\tau>\tau_2$ near the equilibrium points $x_{1,n}^*$. The stable solution near $x_{1,0}^*$ for $\tau=0.4$ is given in Figure \eqref{fig4.7}(b), for $\tau=0.8$ one complex root with positive real part is $0.0322057 + 2.72212 i$ so we get unstable solution and stable solution for $\tau=12$ in Figure \eqref{fig4.7}(d).
				
				\item If $m=8>m_{2*}$ then we are in the SSR region from Figure \eqref{fig4.4}(d). So, we get $\tau_1=0.0173043$ from Figure \eqref{figure4.3}(d) where $g(\tau)=\tau$. So, for $\tau<\tau_1$ we get stable solutions (cf. Figure \eqref{fig4.7}(e) with $\tau=0.01$) and for $\tau>\tau_1$ we get unstable solutions (cf. Figure \eqref{fig4.7}(f) with $\tau=0.03$) near the equilibrium points $x_{1,n}^*$. 
			\end{itemize}
		\end{itemize}	
	\end{Ex}

\begin{figure}
		\subfloat[Stable solution of Example \eqref{eg4.3} for $\tau=1$]{%
			\includegraphics[scale=0.4]{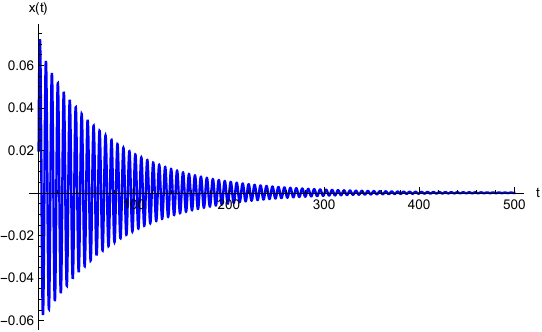}
		}\hspace{0.1cm}
		\subfloat[Unstable solution of Example \eqref{eg4.3} for $\tau=3$]{%
			\includegraphics[scale=0.4]{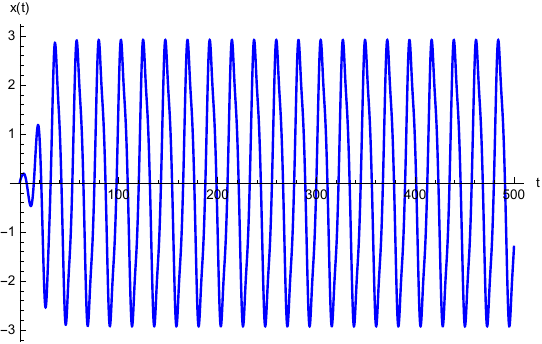}
		}\hspace{0.1cm}
		\subfloat[Untable solution of Example \eqref{eg4.4} for $\tau=2.8$ near the equilibrium point $x_{2,0}^*$]{%
			\includegraphics[scale=0.4]{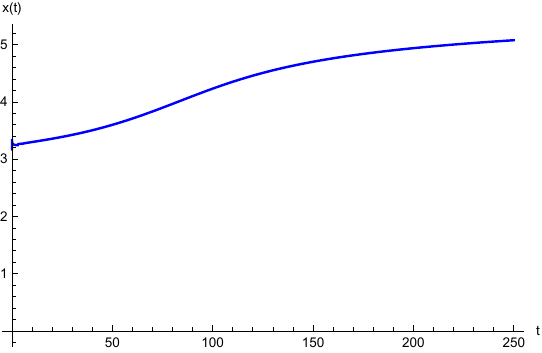}
		}
		\caption{Figures of Example \eqref{eg4.3} and Example \eqref{eg4.4}}
		\label{fig4.9}
	\end{figure}

	\begin{Ex}\label{eg4.3}
		Let us now consider $1/2\leq\alpha<1$. From Subsection \eqref{subsec4.2}, there is no any bifurcation for this case. We get only one critical value $\tau_1$ (see Figure \eqref{figure4.3}(d)) where $g(\tau)=\tau$ such that $0<\tau<\tau_1$ gives stability of $x_{1,n}^*$.
		\begin{itemize}
			\item So, if we fix $\alpha=0.9$, $l=1$ and $m=1.5$  with the initial data $x(t)=0.02$, $\dot{x}(t)=0.1$ for $-\tau<t\leq0$ then we get a critical value $\tau_1=1.03915$ such that for $\tau<\tau_1$ we get stable solution near $x_{1,n}^*$ (cf. Figure \eqref{fig4.9}(a) for $\tau=1$). For $\tau>\tau_1$, we get unstable solution (cf. Figure \eqref{fig4.9}(b) for $\tau=3$).
		\end{itemize}
	\end{Ex}
	
	\begin{Ex}\label{eg4.4}
		The equilibrium points $x_{2,n}^*$ are unstable for all $\tau\geq0$ and $0<\alpha<1$ (Section \eqref{sec4.3}).
			So, if we fix $l=5$, $m=2$, $\tau=2.8$ and $\alpha=0.3$ with the initial condition $x(t)=3.17$,  $(2.8<t\leq0)$ near $x_{2,0}^*=\pi$ then $x_{2,0}^*$ is unstable as shown in Figure \eqref{fig4.9}(c). 	
	\end{Ex}
	
	\section{Conclusions}\label{sec4.7}
	Fractional order generalizations of the classical equations are useful in improving the models. The resulting models are more realistic than their classical counterparts. We generalized the sunflower equation to include two fractional order derivatives. This equation has two sets of countably many equilibrium points viz. $x_{1,n}^*$ and $x_{2,n}^*$. The stability analysis of the equilibrium points is provided by linearizing the equations near the respective equilibrium and using the theory developed in the literature. We proved that the equilibrium points $x_{2,n}^*$ are unstable for all the permissible values of the delay. On the other hand, $x_{1,n}^*$ shows different behaviors described below. For the fractional order $0<\alpha<0.4$, we observed the stable solutions for all the delay parameters and the stability switches for some parameters $l$ and $m$. For $1/2\leq \alpha\leq 1$, we observed a single stable region where the critical value $\tau_1$ of the delay bifurcates the stable behavior from the unstable one.  The fractional order $\alpha=0.4$ shows richer dynamics and we observed all the stability properties and bifurcations described above. The interesting observation is chaos at some parameter sets. The chaotic attractor generated has many scrolls because of the involvement of the sine function in the model.
	
	\section*{Acknowledgment}
	S. Bhalekar acknowledges the University of Hyderabad for the Institute of Eminence-Professional Development Fund (IoE-PDF) by MHRD (F11/9/2019-U3(A)).
	D. Gupta thanks University Grants Commission for financial support (No.F.82-44/2020(SA-III)).
	
	\section*{Declarations} \label{Statements_and_Declarations}
	
	
	\begin{itemize}
		\item {\bf Competing Interests:} The authors declare that there are no competing interests regarding the publication of this research article.
		\item {\bf Data availability:} Data generated in this work is made available at \url{https://drive.google.com/drive/folders/1jOuemmKoSxZfzFSRlotf94YYp5nJI-iy?usp=sharing}.
		\item {\bf Authors' contributions:} S.B. designed the problem. D.G. conducted the research, and both authors collaborated in the discussion. D.G. took the lead in writing the paper, and S.B. provided proofreading and editing. S.B. provided research supervision. 
	\end{itemize}
	\section{Appendix}
	To plot the graphs given in Sections \eqref{sec4.5} and \eqref{sec4.6} we used a numerical method given in the paper \cite{bhalekar2019analysis} obtained from the new predictor and corrector method described in \cite{daftardar2014new}. So, by taking fractional integral of order $2\alpha$ in the equation \eqref{eq4.1} we get, 
	\begin{equation}\label{solving integral equation}
		x(t)=x_0\Big(1+\frac{l t^\alpha}{\tau\Gamma(\alpha+1)}\Big)-\frac{l}{\tau\Gamma(\alpha)}\int_{0}^{t}(t-s)^{\alpha-1}x(s)ds-\frac{m}{\tau\Gamma(2\alpha)}\int_{0}^{t}(t-s)^{2\alpha-1}\sin(x(s-\tau))ds
	\end{equation}
	when $0<\alpha<1/2$ and
	\begin{equation}\label{solving integral equation 1}
		x(t)=x_0\Big(1+\frac{l t^\alpha}{\tau\Gamma(\alpha+1)}\Big)+t x'_0-\frac{l}{\tau\Gamma(\alpha)}\int_{0}^{t}(t-s)^{\alpha-1}x(s)ds-\frac{m}{\tau\Gamma(2\alpha)}\int_{0}^{t}(t-s)^{2\alpha-1}\sin(x(s-\tau))ds
	\end{equation}
	when $1/2\leq\alpha<1$. 
	Note that here $x_0=x(0)$ and $x'_0=x'(0)$.\\
	Now, the next step is to divide the interval $[-\tau,T]$ into $k+N$ sub-intervals where $T/N=\tau/k=h$ the step size. We take a numerical approximation $x_j=x(t_j)$ and for $j\leq0$ we have $x(t_j-\tau)=x(jh-kh)=x_{j-k}$. \\
	Discretizing the equation \eqref{solving integral equation} and using the trapezoidal quadrature formula, we get
	\begin{equation}
		\begin{gathered}
			x_{n+1}=x_0\Big(1+\frac{l t_{n+1}^\alpha}{\tau\Gamma(\alpha+1)}\Big)-\frac{l}{\tau}\frac{h^{\alpha}}{\Gamma(\alpha+2)}\Big(x_{n+1}+\sum_{j=0}^{n} a_{j,n+1} x_j\Big)\\-
			\frac{m}{\tau}\frac{h^{2\alpha}}{\Gamma(2\alpha+2)}\Big(\sin(x_{n+1-k})+\sum_{j=0}^{n} b_{j,n+1} \sin(x_{j-k})\Big)
		\end{gathered}
	\end{equation}
	when $0<\alpha<1/2$ and
	\begin{equation} 
		\begin{gathered}
			x_{n+1}=x_0\Big(1+\frac{l t_{n+1}^\alpha}{\tau\Gamma(\alpha+1)}\Big)+t_{n+1}x'_0-\frac{l}{\tau}\frac{h^{\alpha}}{\Gamma(\alpha+2)}\Big(x_{n+1}+\sum_{j=0}^{n} a_{j,n+1} x_j\Big)\\-
			\frac{m}{\tau}\frac{h^{2\alpha}}{\Gamma(2\alpha+2)}\Big(\sin(x_{n+1-k})+\sum_{j=0}^{n} b_{j,n+1} \sin(x_{j-k})\Big)
		\end{gathered}
	\end{equation}
	when $1/2\leq\alpha<1$.
	Note that
	\begin{equation} 
		a_{j,n+1}=\left\{
		\begin{array}{ll}
			n^{\alpha+1}-(n-\alpha)(n+1)^{\alpha},\quad \textnormal{if}\quad  j=0,\\~\\
			(n-j+2)^{\alpha+1}+(n-j)^{\alpha+1}-2(n-j+1)^{\alpha+1},\quad \textnormal{if}\quad 1\leq j\leq n,\\~\\
			1,\quad \textnormal{if}\quad j=n+1\\~\\
		\end{array}	
		\right\}
	\end{equation}
	and 
	\begin{equation} 
		b_{j,n+1}=\left\{
		\begin{array}{ll}
			n^{2\alpha+1}-(n-2\alpha)(n+1)^{2\alpha},\quad \textnormal{if}\quad  j=0,\\~\\
			(n-j+2)^{2\alpha+1}+(n-j)^{2\alpha+1}-2(n-j+1)^{2\alpha+1},\quad \textnormal{if}\quad 1\leq j\leq n,\\~\\
			1,\quad \textnormal{if}\quad j=n+1.\\~\\
		\end{array}	
		\right\}
	\end{equation}
	So, the predictor terms are
	
	\begin{equation}
		\begin{gathered}
			x_{n+1}^p=x_0\Big(1+\frac{l t_{n+1}^\alpha}{\tau\Gamma(\alpha+1)}\Big)-\frac{l}{\tau}\frac{h^{\alpha}}{\Gamma(\alpha+2)}\Big(\sum_{j=0}^{n} a_{j,n+1} x_j\Big)\\-
			\frac{m}{\tau}\frac{h^{2\alpha}}{\Gamma(2\alpha+2)}\sum_{j=0}^{n} b_{j,n+1} \sin(x_{j-k})
		\end{gathered}
	\end{equation}
	when $0<\alpha<1/2$ and
	\begin{equation}
		\begin{gathered}
			x_{n+1}^p=x_0\Big(1+\frac{l t_{n+1}^\alpha}{\tau\Gamma(\alpha+1)}\Big)+t_{n+1}x'_0-\frac{l}{\tau}\frac{h^{\alpha}}{\Gamma(\alpha+2)}\Big(\sum_{j=0}^{n} a_{j,n+1} x_j\Big)\\-
			\frac{m}{\tau}\frac{h^{2\alpha}}{\Gamma(2\alpha+2)}\sum_{j=0}^{n} b_{j,n+1} \sin(x_{j-k})
		\end{gathered}
	\end{equation}
	 when $1/2\leq\alpha<1$ and
	\begin{equation}
		z_{n+1}^p=\frac{-l}{\tau}\frac{h^\alpha}{\alpha+2}x_{n+1}^p-\frac{m}{\tau}
		\frac{h^{2\alpha}}{\Gamma(2\alpha+2)}\sin(x_{n+1-k}).
	\end{equation}
	The corrector term is
	\begin{equation}
		x_{n+1}^c=x_{n+1}^p-\frac{l}{\tau}\frac{h^\alpha}{\Gamma(\alpha+2)}(x_{n+1}^p+z_{n+1}^p)-\frac{m}{\tau}\frac{h^{2\alpha}}{\Gamma(2\alpha+2)}\sin(x_{n+1-k}).
	\end{equation}
	\bibliographystyle{sn-mathphys-num}
	\bibliography{paperref}

\end{document}